\documentclass[11pt,a4paper]{article}

\usepackage[latin1]{inputenc}

\usepackage{amsfonts,a4wide}
\usepackage{amsmath,mathenv}

\usepackage{epsfig}

\begin{document}
%%%%%%%%%%%%%%
\numberwithin{equation}{section}
\newtheorem{theoreme}{Theorem}[section]
\newtheorem{proposition}{Proposition}[section]
\newtheorem{remarque}{Remark}[section]
\newtheorem{lemme}{Lemma}[section]
\newtheorem{corollaire}{Corollary}[section]
\newtheorem{definition}{Definition}[section]
\renewcommand{\theenumi}{\roman{enumi}}
\def\RR{\mathbb{R} }
%%%%%%%%%%%%%%%%%%%
%%%%%%%%%%%%%%%%%%%%%%%%%%%%%% Pour le fichier de preuve Eric

\def\be{\begin{equation} \displaystyle}
\def\ee{\end{equation} }
\def\ben{\[\displaystyle}
\def\een{\] }
\def\bea{\begin{eqnarray}}
\def\eea{\end{eqnarray} }
\def\bean{\begin{eqnarray*}}
\def\eean{\end{eqnarray*} }
\def\NN{{\rm I\hspace{-0.50ex}N} }
\def\div{{\rm div \;}}
\def\dps{\displaystyle}
\newcommand{\un}{{\mathchoice {\rm 1 \mskip-4mu l} {\rm 1\mskip-4mu
l} {\rm 1\mskip-4.5mu l} {\rm 1\mskip-5mu l}}}

%%%%%%%%%%%%%%%%%%%%%%%%%%%%%%%%%%%%%%%%%%%%%%%%

\title{Well-posedness of a multiscale model\\ for concentrated suspensions}
\author{Eric Canc\`es$^{(a,c)}$,
Isabelle Catto$^{(b)}$,  Yousra
  Gati$^{(a)}$ and Claude Le Bris$^{(a,c)}$ \\
       \footnotesize{(a) CERMICS, Ecole Nationale des
          Ponts et Chauss\'ees,} \\
           \footnotesize{6 \& 8 avenue Blaise Pascal, Cit\'e
          Descartes, 77455 Marne-la-Vall\'ee Cedex 2, France.} \\
\footnotesize{{\tt \{cances,gati,lebris\}@cermics.enpc.fr}} \\
\footnotesize{(b) CEREMADE, UMR CNRS 7534, Universit\'e Paris
 IX-Dauphine,}  \\ \footnotesize{Place du Mar\'echal de Lattre de
 Tassigny, F-75775 Paris
 Cedex 16, France.}\\
\footnotesize{{\tt catto@ceremade.dauphine.fr}} \\
\footnotesize{(c) INRIA Rocquencourt, MICMAC project,
 Domaine de Voluceau,}  \\ \footnotesize{ B.P. 105,
 78153 Le Chesnay Cedex, France.}}

\date{\today}

\maketitle
\noindent{\bf Abstract :}  In a previous work \cite{CCY}, three of us
have studied a nonlinear parabolic equation arising in the mesoscopic
modelling of  concentrated  suspensions of particles that are subjected to a given time-dependent 
shear rate. In the present work we extend the model to allow for a more physically relevant situation when  the shear rate actually depends 
on the macroscopic velocity of the fluid, and as a feedback the
macroscopic velocity 
is influenced by the  average 
stress in  the fluid. The geometry considered is that of a planar Couette
flow. The mathematical system under study couples the one-dimensional
heat equation and a  nonlinear Fokker-Planck type equation with
 nonhomogeneous, nonlocal and possibly
degenerate, coefficients. We show the existence and the 
uniqueness of the global-in-time weak
solution to such a system.

%%%%%%%%%%%%%%%%%%%%%%%%%%%%%%%%%%%%%%%%%%%%%%%%%%%%%%%%%%%%%%
\section{Mechanical context and setting of the equations}

We consider here a concentrated suspension of particles 
in a Couette flow. Examples of such suspensions are numerous: tooth
pastes,  cements, the blood. As opposed to some other complex fluids
such as 
polymeric liquids for which elaborate rheological models, based on fine
mesoscopic physical descriptions, are available (see
e.g. \cite{ottinger}), the modelling of 
concentrated suspensions  is still in its infancy. The specific model
considered 
here however raises interesting mathematical issues, mainly related to
the various nonlinearities present and the coupling of
equations at different scales. Such features are likely to be shared by a
large variety of models, which motivates, and enlarges the scope of,
the present mathematical study.
\vskip6pt
Let us begin with some basics on the mechanical context.
Depending on the concentration,  a suspension of
particles  may 
exhibit different rheological behaviors. At low  concentration, the suspension
behaves like a newtonian fluid at rest or under weak stresses. On the
other hand, when the suspension 
becomes more  concentrated, the motion of each  particle becomes
strongly perturbed by  
the presence of the others and  one observes a so-called jamming
transition where the  
 sample adopts a pastelike behavior. In this transition,  a macroscopic
 yield stress 
 appears~\cite{bi:Larson}. 

  It is well known that when simple fluids are sheared, stress and shear rate are linked
by a linear relation. The linear response coefficients and their
relation to the microstructure of the fluid are well
understood~\cite{hansen}. On the contrary, complex fluids %, such as
                                %concentrated suspensions of hard or
                                %soft spheres,
 exhibit
highly nonlinear  properties far from being understood. These nonlinear properties occur not only at high shear rates, where
one does expect that linear response theory fails, but also at very
low shear rates, which is more surprising.  It is for instance commonly
observed that for some materials (yield stress fluids) the shear stress  
$\sigma$ goes to a non-zero value when the shear rate goes to zero. 

In \cite{bi:HL}, H\'ebraud and Lequeux proposed a model of the
rheological behavior of 
complex fluids based on elementary physical processes.  The system is divided in mesoscopic blocks whose size is large enough
for the stress and strain tensors to be defined for each block. The size
is however small
compared to the characteristic length scale of the stress field. A mesoscopic evolution equation of the stress of each block is then
written:
\begin{enumerate}
   \item at low shear, each particle keeps the same neighbors,
    and a block behaves as an Einstein elastic solid, in which the
    elasticity arises from interactions between neighboring
    particles~;
    \item then, deformation induces local reorganization of
    the particles, at a given stress threshold $\sigma_c$. Above
    this threshold, the block flows as an Eyring fluid~: the
    configuration reached by shearing the suspension relaxes with
    a characteristic time $T_0$ towards a completely relaxed
    state, where no stress is stored~;
    \item lastly, coupling between the flow of neighboring blocks
    must be included. This is taken
    into account by the introduction of a diffusion term in the
    evolution equation, where it is assumed that the diffusion
    coefficient is proportional to the number of reorganizations
    per unit time.
\end{enumerate}
%%%%

%%%%%%

The  equation proposed by H\'ebraud and Lequeux (HL equation in
short) is written for a given shear rate $\dot\gamma$, which only
depends on time:
%%%%%%%
\begin{equation}
  \label{eq:HL-seul}
  \begin{array}{l}
\displaystyle{\partial_t p(t,\sigma)=-G_0\,\dot\gamma (t)
  \;\partial_\sigma p(t,\sigma) \,+\,D(p(t))\,\partial_{\sigma\sigma}^2
p(t,\sigma)} \\ 
\displaystyle{\qquad\qquad\qquad\qquad-\,\frac{\un_{\RR \setminus
    [-\sigma_c,\sigma_c]}(\sigma)}{T_0}\,p 
(t,\sigma)\, +\frac{D(p(t))}{\alpha} \,\delta_0(\sigma)}
  \end{array}
\end{equation}
with 
\begin{equation}
  \label{eq:HL-D}
  D(p(t))=\frac{\alpha}{T_0}\,\int_{|\sigma| > \sigma_c}
  p(t,\sigma) \, d\sigma .
\end{equation}
In the model, each block carries a given shear stress $\sigma$ ($\sigma$ is a real number; it is in fact an extra-diagonal term of the stress tensor in
convenient coordinates).  The evolution of the blocks is described through a probability density $p(t, \sigma)$ which represents the distribution of
stress in the assembly of blocks at time $t$. In equation~(\ref{eq:HL-seul}), $\un_{\RR \setminus [-\sigma_c,\sigma_c]}$
  denotes the characteristic function of the open set $\RR \setminus [-\sigma_c,\sigma_c]$ and $\delta_0$ the Dirac delta function 
  on $\RR$. The three terms arising in the right-hand side of equation
  \eqref{eq:HL-seul} correspond to  the three physical features described above. When a block is submitted to the  shear rate  
$\dot\gamma$, the stress of this block evolves with a variation rate $ G_{0}\,\dot\gamma$  where $G_0$ is an
elasticity constant. When the modulus of the stress  overcomes the critical positive value $\sigma_c$, the block becomes unstable
and may relax into a state with zero stress after a characteristic relaxation time $T_0$. This property is expressed by the last two terms
  in~\eqref{eq:HL-seul}. This relaxation phenomenon induces a rearrangement of the other blocks and this  is finally modelled through the 
  (nonlinear) diffusion term $D(p(t))\;\partial^2_{\sigma\sigma}p$. The diffusion coefficient $D(p(t))$ as given by~\eqref{eq:HL-D}  is assumed to be proportional to
  the density of blocks that relax during time $T_0$. The parameter
  $\alpha$ depends on the microscopic properties of the sample and    is
  supposed to model the ``mechanical  
fragility'' of the material. This nonlinear diffusion term emphasizes
the importance  of collective effects in such materials.
 
\medskip

%%%%%%%
As mentioned above, the shear rate $\dot\gamma$ inserted in
 the original HL equation depends only on time, and not on the space variable. 
 It is
however known from experiment that the shear rate in Couette flows of
non-newtonian fluids is not homogeneous in 
space. In order to better describe the coupling of  the macroscopic flow
 with the evolution of the mesostructure, we therefore introduce a
 space-dependent shear rate $\dot\gamma$ given by the velocity gradient
 (which immediately  implies that an equation of HL type holds
 at \emph{each point} of the sample) and
propose here the following multiscale model for planar Couette flows of concentrated suspensions (see Fig.\ref{fig:Couette} below)~:
\begin{subequations}
 \label{eq:system}
 \begin{EqSystem}
\rho\, \partial_t
U (t,y) = \partial_y\tau (t,y) + \mu\, \partial^2_{yy} U(t,y)\;;\label{eq:u} \\
\nonumber\\
\partial_t p(t,y,\sigma)=-G_0\,\partial_y U (t,y)
  \;\partial_\sigma p(t,y,\sigma) \,+\,D(p(t,y))\,\partial_{\sigma\sigma}^2
p(t,y,\sigma)\,\nonumber \\ 
\qquad\qquad\qquad\qquad-\,\frac{\un_{\RR \setminus
    [-\sigma_c,\sigma_c]}(\sigma)}{T_0}\,p 
(t,y,\sigma)\, +\frac{D(p(t,y))}{\alpha} \,\delta_0(\sigma)\;; \label{eq:p}\\
\tau(t,y)=\int_{\RR}\sigma\,p(t,y,\sigma)\,d\sigma\;;
\label{eq:def-tau}\\
D(p(t,y))=\frac{\alpha}{T_0}\,\int_{|\sigma| > \sigma_c}
  p(t,y,\sigma) \, d\sigma \;;\label{eq:defD}\\
p\,\ge\,0\;;\\
U(0,y)=u_{0}(y)\;, \quad p(0,y,\sigma)=p_0(y,\sigma)\;;\\
 t\in(0;T),\;y\in (0;L),\;\sigma\in \RR\;.
\end{EqSystem}
\end{subequations}
This system is supplied with the initial condition 
 for the probability density:
\begin{equation}\label{IC}
p_0\,\ge\,0\;,\qquad
 \int_{\RR}p_0(y,\sigma)\,d\sigma=1,\quad\text{ for almost every }\; y
\in (0;L)\;.
\end{equation}
and  the
no-slip boundary conditions
\be
U(t,0)=0\,,\quad 
U(t,L)=V(t)\,, \quad \mbox{for almost all } t\;. 
\label{BCu}
\ee
 In the above equations, $U(t,y)$ denotes the component along $e_x$ of the velocity field 
(the flow being laminar and incompressible, the velocity field is of the
 form  $\vec U = U(t,y) \, e_x$), $\rho$ is the volumic mass of the
 fluid and $\mu$ some non-negative viscosity coefficient. 
The function $V$ which appears in the boundary condition~\eqref{BCu}  is a
continuous function on $\RR$ such that $V\in L^\infty(\RR)\cap H^1_{\rm loc}(\RR)$ and $V(0)=0$. 
The initial velocity $u_0$ lies in $L^2(0;L)$.  
\medskip
\begin{figure}[h] \label{fig:Couette}
\centering
\psfig{figure=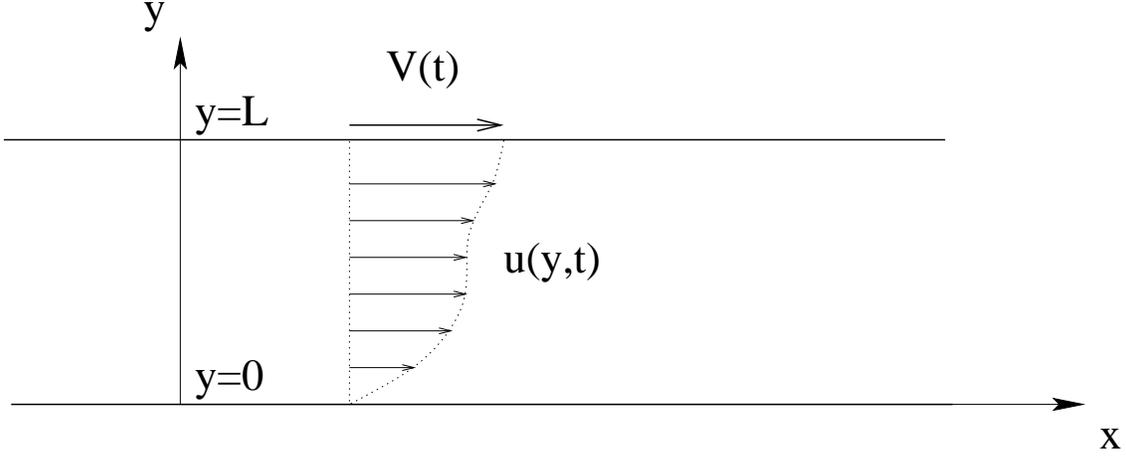,height=6truecm}
\caption{Planar Couette flow}
\end{figure}
\medskip
%%%%%%%%

The mathematical
analysis of the original HL model (\ref{eq:HL-seul}) has been the subject of~\cite{CCY}, 
where a more detailed presentation of the physical background and some additional references may be
found. The main difficulties of course
come from the nonlinearity  in the diffusion term, from the presence of
the singular Dirac mass as a source term,  and foremost from the fact that the parabolic equation 
degenerates if the viscosity coefficient $D(p)$ vanishes. In particular,
it has been shown  that such 
a degeneracy may only occur  if it is already the case at start, that is
$D(p_0)=0$.  And then, the situation becomes very intricate since in
some particular cases several solutions may exist (see~\cite{CCY}). In the coupled
system (\ref{eq:system}), we have to deal with additional difficulties due to the
multiscale coupling: there are \emph{several}
HL-equations (roughly stated, one for each $y$) and all of them are coupled through the
macroscopic equation~\eqref{eq:u}. Proving the well-posedness of the
Cauchy problem for the coupled system is the main purpose of this
article.

\medskip

Before we get to the heart of the matter, we would like to comment on
the diffusion term $\mu\,\partial^2_{yy}u$  in the equation of motion~\eqref{eq:u}. Let us
emphasize that this  artificial viscosity 
has been  added  only for
mathematical   purposes: there is no physical reason why the fluid should be
considered viscous. In the absence of such a regularizing term, we
are  unfortunately unable  to conduct the mathematical analysis in the
whole generality. 

There is
however one particular situation, namely that when $\sigma_c=0$, where we are indeed
able to study the system even without the diffusion term (i.e. with
$\mu=0$). This is the subject of the short Section~3.

\vskip 6pt
\noindent {\bf Notation:}  For given  positive constants $T$ and $L$, we
denote $\Omega=[0;L]$ and $\Omega_T=(0;T)\times\Omega$.  In the sequel,  $C$ a generic positive constant that
may depend on the data but that is independent of $t$, $y$ and
$\sigma$. Also to simplify the notation we shall use in the proofs the
shorthands 
$L^p_T$,  $L^p_y$,  and $L^p_\sigma$  for the
functional spaces $L^p(0;T)$, $L^p(\Omega)$, $L^p(\RR)$, respectively.

%%%%%%%%%%%%%%%%%%%%%%%%%%%%%%%%%%%%%%%%%%%%%%%%%%%%%%%%%%

%%%%%%%%%%%%%%%%%%%%%%%%%%%%%%%%%%%%%%%%%%%%%%%%%%%%%%%%%%
\section{Existence and uniqueness of weak solutions}
\label{general case} 
%%%%%%%%%%%%%%%%%%%%%%%%%%%%%%%%%%%%%%%%%%%%%%%%%%%%%%%% 
In this section, we assume that $\sigma_c>0$, which is the physically
relevant case. With an appropriate change of scales in the coordinates
and variables (see Appendix), we may 
equivalently assume that $L=T_{0}=\sigma_{c}=1$. %  In other words Equations \eqref{eq:u-res} 
% and \eqref{eq:p-res}  are  the ones which are now referred to as \eqref{eq:u} and \eqref{eq:p} respectively. Note that, according to \eqref{eq:u-res}, 
% the coefficient  $\rho$ in \eqref{eq:u}  now represents the Reynolds number. Since the role of the viscosity coefficient $\mu$ in \eqref{eq:u} is only to  
% help for the mathematical analysis, we make the convention that
% $\mu=1$ in the following.
 In addition, without loss of generality, we take $\mu=1$
in~\eqref{eq:u}, for, we recall, $\mu$ is only here for mathematical
convenience and is needed to be strictly positive. 
 
\vskip6pt
%%%%%%%%%%%%%%%%%%%%%%%%%%%%%%%%%%%

%%%%%%%%%%%%%%%%%%%%%%%%%%%%%%%%%%%%

Let us define the velocity field 
$$ \tilde{U}(t,y)\,=\,V(t)\,y\;$$ 
as a lifting
% We look for   a solution $(u;p)$ to our system of
% equations with $\displaystyle{u=u(t,y)}$ in $\mathcal{U_T}$, $\displaystyle{ p=p\,(t,y,\sigma)}$
% in $\displaystyle{L^\infty_{T,y}\big(L^1_\sigma\big)}$
% $\displaystyle{\,\cap\, L^\infty_y(\mathcal{P}_T)}$ and
% $\displaystyle{\sigma\,p}$ in $\displaystyle{
% L^\infty_T\Big(L^2_y\big(L^1_\sigma\big)\Big)}$ (therefore, $\tau\in L^\infty_T(L^2_y)$\/).
% \vskip10pt
% First of all, we shall explain how to construct a lifting
% $\tilde{u} (t,y)$ 
 of the boundary condition (\ref{BCu}), and denote for $T > 0$
$$\mathcal{U}_T=C^0\big([0,T];L^2(\Omega)\big)\cap
L^2\big([0,T];H^1_0(\Omega)\big)$$ 
and
$$\mathcal{P}_T=C^0\big([0,T];L^2(\RR)\big)\cap
L^2\big([0,T];H^1(\RR) \big)\;.$$

% Of course, 
% $\displaystyle{\tilde{u}\,\in\,C^0\big([0,T];L^2(\Omega)\big)\,\cap\,L^2\big([0,  
%  T];H^1_{0}(\Omega)\big)}$.
% (Actually, 
% $\displaystyle{\tilde{u}\,\in\,L^\infty_T\big(C^\infty(\bar{\Omega})\big)}$
% and $\partial_t \tilde{u}\,\in\,L^2_T\times
% C^\infty(\bar{\Omega})$.)

% Moreover,
% \begin{equation*}
% \left\{
% \begin{array}{rll}
% \tilde{u}(0,y)&=&0\;,\quad\quad\;\text{for all}\quad y\,\in[0;1]\;;\\
% \tilde{u}(t,0)&=&0\;,\quad\tilde{u}(t,1)=V(t)\;,\quad\quad\;\,\text{for all}\quad t\,\in[0;T]\;.
% \end{array}
% \right.
% \end{equation*}
% Then, we can recast the system (\ref{eq:u}) by

Setting $\displaystyle{\,U\,=u+\,\tilde{U}}$ and denoting by $\dot V$
the derivative of $V$ with respect to time, the problem
under consideration now reads~:\\

\noindent{\em Find 
  $\displaystyle{u\,\in\, \mathcal{U}_T}$
  and   $\displaystyle{p\,\in\,L^\infty\big([0,T]\times\Omega;L^1(\RR)
  \big)\,\cap\, L^\infty(\Omega;\mathcal{P}_T)}$ solutions
  to
\begin{subequations}\label{eq-u-rel}
  \begin{EqSystem}
\rho\,\partial_t
u\,-\,\partial_{yy}^2u\,=\,\partial_y\tau-\rho\,\dot{V}(t)\,y\;;\label{eq:ubis}\\
\tau=\tau(t,y)=\int_{\RR}\sigma\, p(t,y,\sigma)\, d\sigma\;;\\
 u(0,y)=u_{0}(y)\;;
  \end{EqSystem}
\end{subequations}
\begin{subequations}\label{eq-p-rel}
\begin{EqSystem}
\partial_tp\,+\,G_0\,\Big(\partial_yu\,+\,V(t)\Big)\,\partial_\sigma
p\,-\,D(p(t,y))\,\partial^2_{\sigma\sigma}p\,
+\,\un_{\RR\setminus [-1,1]}(\sigma)\,p\,
=\,\frac{D(p(t,y))}{\alpha}\,\delta_0(\sigma) \;; \label{eq:pbis}
\\ 
p\,\ge\,0\;;\\
D(p(t,y))={\alpha}\,\int_{|\sigma| > 1}
  p(t,y,\sigma) \, d\sigma \;; \\
p(0,y,\sigma)=p_0(y,\sigma)\;.
  \end{EqSystem}
\end{subequations}
}
% These are the equations we shall be considering in the following instead of the
% initial ones \eqref{eq:u} and \eqref{eq:p}.
\vskip10pt\noindent
Our main result is the following~:
%%%%%%%%%%%%%%%%%%%%%%%%%%%%%%%%%%%%%
\begin{theoreme}
\label{main}
%%%%%%%%%%%%%%%%%%%%%%%%%%%%%%%%%%%%%%%
 Let $u_{0}$ be in $L^2(\Omega)$ and let  $p_0$ satisfy the conditions
\begin{equation}\label{ICfull}
\left\{
\begin{array}{l}\dps{
p_0\,\ge\,0\;,\qquad
 \int_{\RR}p_0(y,\sigma)\,d\sigma=1,\quad\text{ for almost every }\; y
\in \Omega\;,}\\
 \dps{p_0\in L^\infty(\Omega\times \RR)\;, \int_\RR
\vert\sigma\vert\,p_0\,d\sigma\;\in\;L^2(\Omega)\;,}
\end{array}
\right. \end{equation} together with:
\begin{equation}\label{strictpos}
\left\{
\begin{array}{l}
\textrm{There exists a positive constant }  \eta\textrm{  such that }\\
\dps{\alpha\,\inf_{\scriptstyle y\in \Omega \atop \scriptstyle \chi\in\RR}
\int_{\vert\sigma+\chi\vert>1}p_0(y,\sigma)\,d\sigma\geq
\eta>0}\;.
\end{array}
\right.
\end{equation}
Then, there exists a unique global-in-time weak solution $(u;p)$ 
$$
u\,\in\,C^0\big(\RR^+;L^2(\Omega)\big)\cap
L^2_{\rm loc}\big(\RR^+;H^1_{0}(\Omega)\big)\;,$$ 
\begin{eqnarray}\label{eq:defP}
    p\in L^\infty\big(\RR^+ \times \Omega ;L^1 (\RR) \big)
\qquad p \in 
L^\infty\big(\Omega;\mathcal{P}_T\big), \quad \forall T > 0 
\end{eqnarray}
to (\ref{eq-u-rel})-(\ref{eq-p-rel}).
In addition, for such a solution,
$p \in L^\infty_{\rm loc} \big(\RR^+ ; L^\infty( \Omega \times\RR)
\big)$ and we have 
\begin{equation*}
\tau\,\in\,L^2\big(\Omega; L^\infty_{\rm loc}(0;T)\big)\cap
C^0\big([0,T];L^2(\Omega)\big)  , 
\end{equation*}
\begin{equation*}
p\,\in\,C^0\big([0,T];L^2(\Omega\times\RR)\big)\quad, \quad\int_\RR p(t,y,\sigma)\,d\sigma=1,
\qquad\mbox{for all $t\geq 0$ and $y\in\,\Omega$\,,}
\end{equation*}
and 
\begin{equation*}
\inf_{\scriptstyle 0\leq t\leq T\atop \scriptstyle y\in \Omega} D(p(t,y))\geq \frac{\eta}{2}\,e^{-T}\;.
\end{equation*}
\end{theoreme}
%%%%%%%%%%%%%%%%%%%%%
Some comments regarding the assumption (\ref{strictpos}) are immediately
in order. 

\medskip

Condition~\eqref{strictpos}  obviously implies that $D(p_0)$ is bounded away from zero independently of $y$.  
The aim of this condition on the initial data $p_0$  is to ensure that, at any time,  the viscosity term
$D(p)$  in~\eqref{eq:pbis} is also  bounded 
away from zero, so that the nonlinear parabolic equation~\eqref{eq:pbis} satisfied by $p$ is non-degenerate at any
time~(see \eqref{eq:nondeg} in Lemma~\ref{lem:borneL2} below). The
condition  is satisfied for example when $p_{0}$ is a Gaussian-like
function. 

Such an assumption seems very demanding, and thus restrictive from the
viewpoint of applications. In fact, some numerical simulations performed
by one of us in~\cite{bi:G} show that even when that assumption
 is not
satisfied at initial time $t=0$, it is indeed satisfied for $t>0$
arbitrarily small. We are unfortunately not able to establish this fact
rigorously, but the numerical evidence mentioned above heuristically shows that
condition  (\ref{strictpos}) can be considered to be always satisfied, up to a
change in the choice of the origin of times.

% \begin{remarque} The subscript $T$ in the functional spaces above means that the norm of the  function
% in the corresponding  space may depend on the length $T$ of the time interval. 
% \end{remarque}

The rest of this section is devoted to the proof of Theorem~\ref{main}.
The  existence and uniqueness  result is first proven on a small time interval with
an argument based on the Banach fixed point Theorem. We introduce the function
${\cal F}_1$ which  associates to every function $u$ in $\dps{L
^2\big([0,T];H^1_{0}(\Omega)\big)}$  the
function $\tau\,=\,\int_\RR \sigma\,p\,d\sigma$ in
$\dps{L^\infty\big([0,T];L^2(\Omega)\big)}$, corresponding to the
(unique) solution $p$ in 
$ L^\infty\big([0,T]\times\Omega;L^1(\RR) \big)
\,\cap\, L^\infty(\Omega;\mathcal{P}_T)$ to \eqref{eq-p-rel}.
Then, we denote by ${\cal F}_2$ the mapping from
$\dps{L^\infty\big([0,T];L^2(\Omega)\big)}$ to
$\dps{L^2\big([0,T];H^1_{0}(\Omega)\big)}$, which 
associates to every $\tau$ in
$\dps{L^\infty\big([0,T];L^2(\Omega)\big)}$, the unique 
solution $v$ (in $\mathcal{U}_T$\,) to the heat equation~:
\begin{subequations}\label{eq:eqv}
 \begin{EqSystem}
\rho\,\partial_t v-\,\partial^2_{yy}v\,=\,\partial_y\tau-\rho\,\dot{V}(t)\,y\;\quad\text{on
}\;\Omega_T\;;\\ v(0,y)=u_0 \; .
 \end{EqSystem}
\end{subequations}
We next define the mapping  ${\cal F}$ on
$\dps{L^2\big([0,T];H^1_{0}(\Omega)\big)}$  as ${\cal F}={\cal F}_2\circ\,{\cal F}_1$~:
\begin{equation}
\begin{array}{lcccccc}
&{\cal F}:&L^2\big([0,T];H^1_{0}(\Omega)\big)&\buildrel{{\cal
F}_1}\over\longrightarrow &L^\infty\big([0,T];L^2(\Omega)\big) &\buildrel{{\cal
F}_2}\over\longrightarrow& L^2\big([0,T];H^1_{0}(\Omega)\big)\\
& &u&\longmapsto&\tau&\longmapsto &v
\end{array}
\end{equation}
and our main step consists in  proving  the following. \vskip6pt
%%%%%%%%%%%%%%%%%%%%%%%%%%%%%%%%%%%%%%%%%%%%%%%%%%%%%%%%%%%%%%%%%%%%%%
\begin{proposition}
\label{prop:ptfixe} For every $T>0$, the mapping ${\cal F}$ from
$\dps{L^2\big([0,T];H^1_{0}(\Omega)\big)}$ into itself is  well-defined and for
$T>0$ small enough,  it admits a unique fixed point denoted by
$u$.
\end{proposition}
%%%%%%%%%%%%%%%%%%%%%%%%%%%%%%%%%%%%%%%%%%%%%%%%%%%%%%%%%%%%%%%%%%%%%%%

\vskip10pt

\noindent The proof  of  Proposition~\ref{prop:ptfixe} is organized as follows.
We first check in Lemma~\ref{lem:F2} below that ${\cal F}_2$ is well-defined and that it is 
a Lipschitz continuous function with a Lipschitz constant that may be
chosen arbitrarily small provided the length of the time interval is reduced~(Section \ref{sec:f2_contraction}).  In Section~\ref{sec:f1_bien_def}, we next prove that ${\cal F}_1$ is well-defined. 
We establish  in Section~\ref{sec:f1_lipschitz}  that ${\cal F}_1$ is a Lipschitz continuous function with a locally bounded Lipschitz constant 
with respect to time interval. Therefore the composed mapping
$\mathcal{F}$  is contracting on small enough time interval. The  
existence and uniqueness of a solution on a small time interval follows  by the Banach fixed point theorem.  Finally in
Section~\ref{sec:longtime} we deduce the existence and uniqueness of the
global-in-time solution.

\medskip

Henceforth, and unless otherwise stated, the initial condition $p_0$ is
fixed and it satisfies the assumptions \eqref{ICfull}-\eqref{strictpos}
of the statement of the Theorem. 

%%%%%%%%%%%%%%%%%%%%%%%%%%%%%%%%%%%%%%%%%%%%%%%%%%%%%%%%%%%%%%%%%
\subsection{The map ${\cal F}_2$  is  a contraction on $[0,T]$ for $T$
  small enough}
\label{sec:f2_contraction}
%%%%%%%%%%%%%%%%%%%%%%%%%%%%%%%%%%%%%%%%%%%%%%%%%%%%%%%%%%%%%%%%%%%%%%%%%%%%%
%%%%%%%%%%%%%%%%%%%%%%%%%%%%%%%%%%%%%%%%%%%%%%%%%%%%%%%%%%%%%
\begin{lemme}
\label{lem:F2}
%%%%%%%%%%%%%%%%%%%%%%%%%%%%%%%%%%%%%%%%%%%%%%%%%%%%%%%%%%%%%%%%%%%%%%%%%
For every $T>0$, the mapping ${\cal F}_2$ is  Lipschitz continuous
from $\dps{L^\infty\big([0,T];L^2(\Omega)\big)}$ to
$\dps{L^2\big([0,T];H^1_{0}(\Omega)\big)\,} $, and the Lipschitz constant
goes to $0$  with $T$.
\end{lemme}
\vskip6pt\noindent\textbf{Proof of Lemma \ref{lem:F2}:} 
We first observe that the mapping ${\cal F}_2$ is well-defined.
Indeed, for every function $\tau$ in $L^\infty([0,T];L^2(\Omega))\,$,\,
$\partial_y\tau\,\in\,L^\infty([0,T];H^{-1}(\Omega))$, and therefore, the
existence and the uniqueness of a solution
$v\,\in\,\mathcal{U}_T$ of the
heat equation (\ref{eq:eqv}) is a  standard result. Let now
$\tau_1$ and $\tau_2$ be two functions in $L^\infty([0,T];L^2(\Omega))\,$,
and let us denote
 $v_1={\cal F}_2(\tau_{1})$ and
 $v_2={\cal F}_2(\tau_{2})$. We also
 set   $v\,=\,v_1-v_2$ and
$\tau\,=\,\tau_1-\tau_2$. Then, $v$ satisfies
\begin{subequations}
\label{eq:v}
\begin{EqSystem}
\rho\,\partial_t v-\,\partial^2_{yy}v\,=\,\partial_y\tau\quad\text{on }\;\Omega_T\;;\\
v(0,y)=0\;;\\
v(t,0)=v(t,1)=0\;,
 \end{EqSystem}
 \end{subequations}
and if we apply Equation \eqref{eq:v} to $v$ and integrate
over $\Omega$ we get
\begin{equation}\label{eq:v2}
\frac{\rho} 2 \frac{d}{dt}\int_\Omega \vert v\vert^2\,+\,\int_\Omega \vert
\partial_y v\vert^2\,=\,-\int_\Omega \tau \,\partial_y v\;.
\end{equation}
By the Cauchy-Schwarz  and the Young inequalities, we obtain for $t\,\in\,[0;T]$,
$$\rho\,\int_\Omega \vert v\vert^2\,+\, \int_0^t\Big(\int_\Omega \vert
\partial_y v\vert^2\,dy\Big)\,ds\,\le\, \int_0^T\Big(\int_\Omega \vert\tau\vert^2\,dy\Big)\,ds\;,$$
and therefore by the Poincar\'e inequality
\begin{equation}\label{eq:lipF2}
\Vert
v\Vert_{L^2([0,T];H^1(\Omega))}\, \le 
2 \sqrt{T} \Vert \tau \Vert_{L^\infty([0,T];L^2(\Omega))}\;. 
\end{equation}
\vskip6pt
\null\hfill$\diamondsuit$\\

%%%%%%%%%%%%%%%%%%%%%%%%%%%%%%%%%%%%%%%%%%%%%%%%%%%%%%%%%%%%%%%%%%%%
\subsection{The map  ${\cal F}_1$ is well defined}
\label{sec:f1_bien_def}
%%%%%%%%%%%%%%%%%%%%%%%%%%%%%%%%%%%%%%%%%%%%%%%%%%%%%%%%%%%%%%%%%%%%
\vskip10pt

\noindent Equation~(\ref{eq-p-rel})
with the variable $y$ frozen has been studied in \cite{CCY}. For the
sake of consistency we now recall~:
%%%%%%%%%%%%%%%%%%%%%%%%%%%%%%%%%%
\begin{proposition}\cite[Theorem 1.1]{CCY} {\bf (Global-in-time
    existence for all $y$)}
\label{prop:main1}
%%%%%%%%%%%%%%%%%%%%%%%%%%%%%%%%%%

For
almost every $y$ in $\Omega$, let $b(\cdot ,y)$ be a given function in
$L^2_{\rm loc}(\RR^+)$, and let $p_0$ such that~:
\begin{equation}\label{eq:IC-full}
p_0(y,\cdot)\in L^1(\RR) \cap L^\infty(\RR) \,,\quad p_0(y,\cdot)\geq 0\,,\quad \int_\RR
p_0(y,\sigma)\,d\sigma=1\quad\mathrm{ and }\int_\RR|\sigma|\,p_0\,d\sigma<+\infty\;,
\end{equation}
and
\begin{equation}
\label{eq:p0-prop}
D(p_0(y))>0\;.
\end{equation}
Then, for every $T>0$ and for almost every $y$ in
$\Omega$, there exists a unique solution $p=p(t,y,\sigma)$ in
$L^\infty([0,T];L^1(\RR) \cap L^\infty(\RR))\cap L^2([0,T];H^1(\RR))$ to
the equation
\begin{subequations}\label{eq:syst-p}
 \begin{EqSystem}
    \partial_{t}p=-b(t,y) \; \partial_{\sigma}p+D(p(t,y))\,
    \partial^2_{\sigma\sigma}p-\,\un_{\RR\setminus
    [-1,1]}(\sigma)\,p+\frac{D(p(t,y))}{\alpha}\;\delta_{0}(\sigma) \;; \quad 
\\
p\,\ge\,0\;;\\D(p(t,y))={\alpha}\,\int_{|\sigma| > 1}
  p(t,y,\sigma) \, d\sigma \;; \\
p(0,y, \sigma)=p_0(y,\sigma)\;.
  \end{EqSystem}
\end{subequations}
In addition, for almost every $y$ in $\Omega$, we have
\begin{itemize}
\item 
$\dps \int_\RR p(t,y, \sigma) \, d\sigma=1\qquad \textrm{ for all }\, t
\geq 0\;,$
\item for all $T > 0$,
\begin{equation}\label{eq:bd-infty}
\max_{0\leq t\leq T}\Vert p(t,y, \cdot) \Vert_{L^\infty_{\sigma}}\leq  \Vert
p_0(y,\cdot)\Vert_{L^\infty_\sigma}+\frac{\sqrt{\alpha}\,\sqrt{T}}{\sqrt\pi}\;, 
\end{equation}
\item 
$p(\cdot,y)\in C^0([0,T];L^1(\RR) \cap L^2(\RR))$,
\item 
 $D(p(\cdot,y))\in C^0([0,T])$, 
\item 
 for
every $T>0$ there exists a positive constant $\eta(T,y)$ such that
\begin{equation}\label{eq:nondeg-weak}
\min_{0\leq t\leq T}D(p(t,y))\geq \eta(T,y)\;.
\end{equation}
\item for almost all $y$, 
  $(t,\sigma) \mapsto \sigma\,p( t, y,\sigma ) \in
  L^\infty([0,T];L^1(\RR))$, so that  
the average stress $\tau(\cdot,y)$ is well defined by (\ref{eq:def-tau})
  in $L^\infty_{\rm loc}(\RR^+)$.\end{itemize} 
\end{proposition}
\medskip
\noindent We now fix some initial condition $p_0$ satisfying the
conditions~\eqref{ICfull} and \eqref{strictpos} (thus \emph{a fortiori}
the conditions (\ref{eq:IC-full}) and (\ref{eq:p0-prop})) and set 
$$b(t,y)=G_0\,(\partial_y u(t,y)+V(t))
$$
for $u\in L^2([0,T];H^1_0(\Omega))$. In view of Proposition~\ref{prop:main1}, we
know the  existence and uniqueness of a solution $p$ to~\eqref{eq:p} for
given $u$. Our next step now consists in analyzing the dependence
on the variable $y$.

%%%%%%%%%%%%%%%%%%%%%%%%%%%%%%%%%%%%%%%%%%%%%%%%%%%%%%%%%%%%%%%%%%%%
\begin{lemme}{\bf (Uniform-in-$y$ \emph{a priori} estimates on $p$)}
\label {lem:borneL2}
%%%%%%%%%%%%%%%%%%%%%%%%%%%%%%%%%%%%%%%%%%%%%%%%%%%%%%%%%%%%%%%%%%%%
Let $T>0$ be given. We assume that the initial data $p_0$
satisfies \eqref{ICfull} and
\begin{equation*}
\inf_{y\in\Omega}D(p_{0}(y))>0\;.
\end{equation*}
Notice that (\ref{strictpos}) is not needed, but (\ref{eq:IC-full})
and (\ref{eq:p0-prop})) are fulfilled. Then, if we denote by $p$ the unique solution to \eqref{eq:syst-p}
given by Proposition~\ref{prop:main1}, we have~:
\begin{enumerate}
\item[(i)] $p\,\in\,L^{\infty}([0,T]\times\Omega;L^1_{\sigma}\cap
  L^\infty_{\sigma})$ with 
\begin{equation}\label{eq:borneLinfty}
\Vert
p\Vert_{L^\infty([0,T]\times\Omega\times\RR)}\,\le\,\Vert
p_0\Vert_{L^\infty(\Omega\times\RR)}+\frac{\sqrt\alpha\;\sqrt{T}}{\sqrt\pi}
\end{equation}
and 
\begin{equation}\label{eq:intp}
\int_\RR p(t,y, \sigma) \, d\sigma=1\qquad \textrm{ for all }\, t
\geq 0, \,\textrm{ for almost every } y\textrm{ in } \Omega\;. \end{equation}
%  therefore
% \begin{equation}\label{eq:borneL2}
% \sup_{\scriptstyle 0\leq t\leq T\atop \scriptstyle y\in\Omega}\,\Vert
% p\Vert_{L^2_{\sigma}}^2\,\le\,\Vert
% p_0\Vert_{L^\infty_{y,\sigma}}+\frac{\sqrt\alpha\;\sqrt{T}}{\sqrt\pi}\;.
% \end{equation}
\item[(ii)] The stress $\tau$ is in $L^2(\Omega,L^\infty([0,T]))$ (hence in $L^\infty([0,T];L^2(\Omega))$\/). 
\item[(iii)] If in addition $p_{0}$ satisfies the non-degeneracy condition \eqref{strictpos}, we have
\begin{equation}\label{eq:nondeg}
\inf_{\scriptstyle0\leq t\leq T\atop\scriptstyle y\in \Omega}\,D(p(t,y))\geq \frac{1}{2}\,e^{-T}\,\eta\;,
\end{equation}
and
 \begin{equation}\label{eq:bdgrad}
\sup_{y\in\Omega}\int_{0}^T\Big(\int_{\RR}\vert\partial_{\sigma}p\vert^2\,d\sigma\,\Big)\,dt
\leq \frac{2}{\eta}\,e^{T}\,\left(\Vert
p_0\Vert_{L^\infty(\Omega\times\RR)}\Big(\frac{1}{2}+T\Big)+
\frac{\alpha\;}{\sqrt\pi}\,T^{3/2} \right)\;.
\end{equation}
\end{enumerate}
\end{lemme}
\vskip6pt

\noindent\textbf{Proof of Lemma \ref{lem:borneL2}:} 

\noindent To prove Assertion (i), we use
the estimates obtained in~\cite{CCY} with  the variable $y$ kept frozen. The assumptions on $p_0$ ensure  that $p$ is in
$L^\infty_{T,\sigma}$ for almost every $y$, that \eqref{eq:intp} holds,  and that \eqref{eq:bd-infty} holds by virtue of \cite[Proposition 1.1, Eq.(1.8)]{CCY}. 
Estimate \eqref{eq:borneLinfty} follows. % By virtue of~\eqref{eq:intp} in Proposition~\ref{prop:main1}, \eqref{eq:borneL2} immediately follows by interpolation. 

Assertion (ii) follows from \cite[Proposition 1.1, Eq.(1.9)]{CCY}~: for almost every $y$ in $\Omega$,
\begin{eqnarray*}
\sup_{0\leq t\leq T}\int_\RR\vert\sigma\vert\, p\,d\sigma\leq
\int_\RR \vert\sigma\vert\,
p_0\,d\sigma+\sqrt{T}\,\Big(\frac{2\,\sqrt{\alpha}}{\sqrt{\pi}}+\Vert
b\Vert_{L^2_T}\Big)+\frac{2}{3}\,T^{3/2}\,\Big(1+\frac{2\sqrt{\alpha}}{\sqrt\pi}\Big)\;,
\end{eqnarray*}
with  $b=b(t,y)=G_{0}\,\big(\partial_yu(t,y) +V(t)\big)$. Then
\begin{equation*}\label{eq:borne_sur_sp} 
\Big\Vert\sup_{0\le t\le
T}\int_\RR\vert\sigma\vert
p\,d\sigma\;\Big\Vert_{L^2_y}\,\le\,\Vert \sigma\,
p_0\Vert_{L^2_y(L^1_\sigma)}\,+\,C(T, \alpha,\Vert
V\Vert_{L^\infty([0,T])})\,+\,\sqrt{T}\,G_{0}\,\Vert\partial_y
u\Vert_{L^2(\Omega_T)}\;,
\end{equation*}
with 
\begin{equation*}
C(T,\alpha,\Vert V\Vert_{L^\infty([0,T])})=G_{0}\,T\,\Vert
V\Vert_{L^\infty([0,T])}+\frac{\sqrt{T}}{\sqrt\pi}
\,\Big(2\,\sqrt{\alpha}+\frac{2}{3}\,T \,\big(\sqrt \pi+
2\sqrt{\alpha}\big)\Big)\;.
\end{equation*}

 For Assertion (iii), following~\cite[Proof of Lemma 3.1]{CCY}, 
we define 
\begin{equation*}%\label{def:t*}
t^*(y) =\inf\left\{t>0\,; \int_{\big\vert\sigma+\int_0^t b(s,y)\,ds\big\vert>1} p_0(y,\sigma)\,d\sigma=0\,\right\}\;,
\end{equation*}
Because  $p_{0}$ satisfies the non-degeneracy
condition~\eqref{strictpos}, we have for all $y$,  $t^*(y)>0$  and  the support of $p_0$
is contained in the interval 
$\left[-1-\int_0^{t^*(y)} b(s,y)\,ds,1-\int_0^{t^*(y)}
  b(s,y)\,ds\right[$.

 Moreover, 
 for any $T=T(y)<\frac{t^*(y)}{2}$,
\begin{equation}\label{eq:nu1}
\min_{0\leq t\leq T}D(p(t,y))\geq \frac{\alpha}{2}\,e^{-T}\,\min_{0\leq t\leq T}
 \,\int_{\left\vert\sigma+\int_0^{t} b(s,y)\,ds\right\vert>1} p_0(y,\sigma)\,d\sigma\;.
\end{equation}
% (The proof of this fact is actually reproduced below in the course of
% the proof of Lemma~\ref{lem:pnondeg}.) 

The assumption  \eqref{strictpos} on $p_0$  ensures that   $t^*(y)=+\infty$ for almost every $y$ in $\Omega$. 
Therefore \eqref{eq:nu1} holds true on any time interval $T>0$ independently of $y$ and \eqref{eq:nondeg} is an immediate consequence of \eqref{eq:nu1} by using \eqref{strictpos}. 
\vskip6pt
Finally \eqref{eq:bdgrad} follows in a very standard way  from
\eqref{eq:nondeg} and \cite[Equation (3.7)]{CCY}, multiplying  \eqref{eq:p} by $p$,  integrating  over $\RR$ with respect to $\sigma$, and   using  
\eqref{eq:nondeg} and the previous bounds. \hfill$\diamondsuit$
\noindent
%%%%%%%%%%%%%%%%%%%%%%%%%%%%%%%%%%%%%%%%%%%%%%%%%%%%%%%%%%%%%%%%%%%%%%%%%%%%%
\subsection{The map ${\cal F}_{1}$ is Lipschitz continuous }
\label{sec:f1_lipschitz}
%%%%%%%%%%%%%%%%%%%%%%%%%%%%%%%%%%%%%%%%%%%%%%%%%%%%%%%%%%%%%%%%%%%%%%%%%%%%%%
%%%%%%%%%%%%%%%%%%%%%%%%%%%%%%%%%%%%%%%%%%%%%%%%%%%%%%%%%%%%%
\begin{lemme}
\label{lem:F1}
%%%%%%%%%%%%%%%%%%%%%%%%%%%%%%%%%%%%%%%%%%%%%%%%%%%%%%%%%%%%%%%%%%%%%%%%%
For every $T>0$, the mapping ${\cal F}_1$ is  Lipschitz continuous
from $\dps{L^2\big([0,T];H^1_{0}(\Omega)\big)\,} $ to
$\dps{L^\infty\big([0,T];L^2(\Omega)\big)}$, and the Lipschitz constant
is a locally bounded function of $T$.
\end{lemme}
\vskip6pt\noindent\textbf{Proof of Lemma \ref{lem:F1}:} Let $u_1$
and $u_2$ be two functions in $\dps{L^2\big([0,T];H^1_{0}(\Omega)\big)} $, and let  
$\tau_1={\cal F}_1(u_{1})$ and $\tau_2={\cal F}_1(u_{2})$. We denote by
$p_i$, $i=1,\,2$, the unique  solution to \eqref{eq-p-rel}  corresponding to
$u_i$ whose existence is guaranteed by Proposition~\ref{prop:main1} 
and Lemma~\ref{lem:borneL2}.  We also set $v=u_1-u_2$,  $q\,=\,p_1-p_2$ and
$\tau\,=\,\tau_1-\tau_2$. Recall that, for $i=1,\,2$, $\tau_i\,=\,\int_\RR \sigma\,p_i\,d\sigma$. We formally multiply equation (\ref{eq-p-rel}) by $\sigma$ and
integrate it over $\RR$ with respect to $\sigma$ to find
\begin{equation*}
\left\{ \begin{array}{l}
\dps{\partial_t \tau_i\,+\,\tau_i\,=\,G_0\,\big( \partial_y
u_i+V(t)\big)+\,\int_{|\sigma|\le 1} \sigma \,p_i\,d\sigma }\;;\\
\dps{\tau_i(0,y)\,=\,\int_{\RR}\sigma\,p_{0}\,d\sigma\;.}
\end{array} \right.
\end{equation*}
The argument may be made rigorous with the help of a standard cut-off argument as in ~\cite{CCY}. Subtracting the equations satisfied by $\tau_1$ and $\tau_2$
yields
\begin{equation}\label{eq:entau}
\left\{ 
\begin{array}{l}
\dps{\partial_t \tau\,+\,\tau\,=\, G_0\,\partial_y
v\,+\,\int_{|\sigma|\le 1} \sigma\, q\,d\sigma \;; }\\
\tau(0,y)\,=\,0\;.
\end{array} 
\right.
\end{equation}
We then apply $\tau$ to  (\ref{eq:entau}) and integrate over
$\Omega$ to obtain 
\begin{equation}\label{eq:tau2}
\frac 1 2 \frac{d}{dt} \int_\Omega|\tau|^2\,+ \,\int_\Omega
|\tau|^2\,=\,G_0\,\int_\Omega \partial_y
v\;\tau\,+\,\int_\Omega\tau\,\left(\int_{|\sigma|\le 1} \sigma \,q\,d\sigma\right)\,dy \;.
\end{equation}
Using the Young and the Cauchy-Schwarz  inequalities, we have
\begin{equation*}
\int_\Omega |\partial_y v\,\tau|\leq \frac{1}{2\,G_0}\Vert\tau\Vert^2_{L^2(\Omega)}\,+\,\frac{G_0}{2}\Vert\partial_y
v\Vert^2_{L^2(\Omega)}
\end{equation*}
and
\begin{eqnarray}
\left|\int_\Omega\tau\,\left(\int_{|\sigma|\le 1} \sigma \,q\,d\sigma\right)\,dy\right|
&\le&\int_\Omega \;\Big(\int_\RR\vert
 q\vert^2\,d\sigma\Big)^{1/2}\, |\tau|\,\,dy\nonumber\\
&\le&\frac{1}{2} \int_\Omega\Big(\int_\RR\vert q\vert^2\,d\sigma\Big)\,dy\,+\,\frac{1}{2}\Vert\tau\Vert^2_{L^2(\Omega)}
\;,\label{eq:evalrhs}
\end{eqnarray}
thus
\begin{equation*}
\frac{d}{dt} \int_\Omega|\tau|^2\,\le\,G_0^{\,2}\,\Vert\partial_y
v\Vert^2_{L^2(\Omega)}\,+\,\int_\Omega\Big(\int_\RR\vert q\vert^2\,d\sigma\Big)\,dy\;,
\end{equation*}
and
\begin{equation}
\label{eq:ineq-tau}
\sup_{0\leq t\leq T}\Vert \tau\Vert^2_{L^2(\Omega)}\,\le\,G_0^{\,2}\,\Vert\partial_y
v\Vert^2_{L^2(\Omega_T)}\,+\,\Vert q \Vert^2_{L^2(\Omega_T\times\RR)}\;.
\end{equation}
Let us now admit for a while that
\begin{equation}
\label{eq:estiq}
\Vert q \Vert^2_{L^2(\Omega_T\times\RR)}
\,\le\,C(T)\, \Vert
\partial_y\,v\Vert^2_{L^2(\Omega_T)}\;,
\end{equation}
with  $C(T)$ being a locally bounded function of $T$.  Inserting
\eqref{eq:estiq} into (\ref{eq:ineq-tau}), we obtain
\begin{equation}
\label{cstlip}
\sup_{0\leq t\leq T}\Vert \tau\Vert^2_{L^2(\Omega)}\,\le\,\big(G_0^{\,2}\,+\,C(T)\big)    \Vert\partial_y
v\Vert^2_{L^2(\Omega_T)}\;,
\end{equation}
and therefore, the mapping ${\cal F}_1$ is indeed Lipschitz continuous.\\

In order to establish \eqref{eq:estiq}, we  subtract the
equations (\ref{eq-p-rel}) satisfied by $p_{1}$ and $p_{2}$ respectively to deduce that $q$ solves   
\begin{subequations}\label{eq:q}
\begin{EqSystem}
\partial_t q=-\,G_0\,\partial_y u_1\,\partial_{\sigma}q\,-\,G_0\,\partial_y  v\,\partial_{\sigma}p_2\,
-\,G_0\,V(t)\,\partial_{\sigma}q+D(q)\,\partial_{\sigma\sigma}^2  p_1\nonumber\\
\qquad\qquad+D(p_2)\,\partial_{\sigma \sigma}^2 q
-\un_{\RR\setminus [-1,1]}(\sigma)\,q  +\,\frac{D(q)}{\alpha}\,
\delta
      _0(\sigma)\;;\\
q(0,y,\sigma)=0\;,
\end{EqSystem}
\end{subequations}
for almost every $y$ in $\Omega$. Then, we apply (\ref{eq:q}) to $q$, and 
integrate with respect to $\sigma$, to obtain
\begin{eqnarray}
 \lefteqn{\frac{1}{2}\,{\partial \over{\partial t}}\int_{\mathbb{R}} q 
^2+\int_{\vert\sigma\vert>1}q^2\,
+D(p_2)\,\int_{\mathbb{R}}|\partial_{\sigma}q|^2\,}\nonumber\\
&=& \,G_0\,\partial_y v\int_{\mathbb{R}} p_2\,\partial_{\sigma}
 q-D(q)\int_{\mathbb{R}}\partial_{\sigma}p_1\;\partial_{\sigma}q \,+\,\frac{D(q)}{\alpha}\,q(t,y,0)\;.\label{M1}
\end{eqnarray}
By the Cauchy-Schwarz and the Young inequalities and  using the bound from below~\eqref{eq:nondeg} on $D(p_{2})$,  we have
\begin{eqnarray}
\left\vert\partial_y v\int_{\RR} p_2\,\partial_{\sigma}
 q\right\vert&\le&\vert\partial_y v \vert\;\Vert p_2\Vert
_{L^2_{\sigma}}\;\Vert
 \partial_{\sigma}q \Vert_{L^2_{\sigma}}\nonumber\\
&\le&\frac{3\,G_0}{4\,D(p_2)}\Vert p_2\Vert
_{L^2_{\sigma}}^2\;\vert\partial_y v
 \vert^2+\,\frac{D(p_2)}{3\,G_0}\,\Vert\partial_{\sigma}q
 \Vert_{L^2_{\sigma}}^2\nonumber\\
&\leq& \frac{3\,G_0}{2\,\eta}\;e^{T}\;\Big(\Vert p_{0}\Vert_{L^\infty_{y,\sigma}}
+\frac{\sqrt{\alpha\,T}}{\sqrt{\pi}}\,\Big)\vert\partial_y v
 \vert^2+\,\frac{D(p_2)}{3\,G_0}\,\Vert\partial_{\sigma}q
 \Vert_{L^2_{\sigma}}^2\;, \label{term1}
\end{eqnarray}
thanks to the $L^\infty$ bound \eqref{eq:borneLinfty} on $p_2$ and the fact
that $\dps \int_{\RR} p_2(t,y,\sigma) \, d\sigma=1$.  In a similar way,
we obtain 
\begin{eqnarray}
\vert  D(q)\int_{\RR}\partial_{\sigma}p_1\,\partial_{\sigma}q\vert
&\le&\,\vert D(q) \vert\; \Vert
\partial_{\sigma}p_1 \Vert _{L^2_{\sigma}}\;\Vert
\partial_{\sigma}q
\Vert_{L_{\sigma}^2}\nonumber\\
&\le&  \frac{3}{2\,\eta}\;e^{T}\;\vert D(q) \vert^2\; \Vert
\partial_{\sigma}p_1 \Vert _{L^2_{\sigma}}^2+\,\frac{D(p_2)}{3}\;\Vert
\partial_{\sigma}q
\Vert_{L_{\sigma}^2}^2\;.\label{Y1}
\end{eqnarray}
As  $\dps{\int_{\RR }q=0}$ (recall $\dps{\int_{\RR}p_{1}=\int_{\RR}p_{2}=1}$),  we may write
\begin{equation}
\vert D(q) \vert=\alpha\Big\vert\int_{|\sigma|\geq 1}\,q\Big\vert
=\alpha \Big\vert\,\int_{|\sigma|\le1}\,q\,\Big\vert\le\alpha\,\sqrt{2}\;
\Vert\,q\,\Vert_{L^2_\sigma}\;.\label{eq:N}
\end{equation}
Thus, inserting (\ref{eq:N})  into \eqref{Y1}, we obtain
\begin{equation} \label{term2}
\vert D(q)\int_{\RR}\partial_{\sigma}p_1\,\partial_{\sigma}q\vert
\,\le\,\frac{3\,\alpha^{2}}{\eta}\;e^{T}\;\Vert\,q\,\Vert^2_{L^2_\sigma}\,\Vert\,\partial_\sigma
p_1\,\Vert^2_{L^2_\sigma}\,+\,\frac{D(p_2)}{3}\;\Vert\,
\partial_{\sigma}q\,\Vert_{L_{\sigma}^2}^2.
\end{equation}
On the other hand, from  the Sobolev embedding
$H^1(\RR)\,\hookrightarrow\,L^\infty(\RR)$, we know
 \begin{equation}\label{eq:sob}
 \Vert q \Vert_{L^\infty_{\sigma}}\leq\frac{1}{\sqrt 2}\, (\parallel
q\parallel_{L^{2}_{\sigma}}^2+\Vert\partial_{\sigma}q
\Vert_{L^2_{\sigma}}^2)^{\frac{1}{2}}\;,
\end{equation}
and next using successively  Young's inequality, (\ref{eq:N}),
\eqref{eq:nondeg} again and the fact that $D(p_2)\leq\alpha$, we find
\begin{eqnarray}
\frac{1}\alpha\,\vert D(q) \vert\,\vert
q(t,y,0) \vert&\le&\Vert\,q\,\Vert_{L^2_\sigma}(\parallel
q\parallel_{L^{2}_{\sigma}}^2+\Vert\partial_{\sigma}q
\Vert_{L^2_{\sigma}}^2)^{\frac{1}{2}}\nonumber\\
&\le&\Big(\frac{D(p_{2})}{3}+\frac{3}{4\, D(p_{2})}\Big)\;\Vert q
\parallel_{L^{2}_{\sigma}}^2\,+\,\frac{D(p_2)}{3}\,\Vert\partial_{\sigma}q
\Vert_{L^2_{\sigma}}^2\nonumber\\
&\le&\Big(\frac{\alpha}{3}+\frac{3}{2\,\eta}\,e^{T}\Big)\;\Vert q
\parallel_{L^{2}_{\sigma}}^2\,+\,\frac{D(p_2)}{3}\,\Vert\partial_{\sigma}q
\Vert_{L^2_{\sigma}}^2\;.\label{term3}
\end{eqnarray}
Inserting \eqref{term1}, \eqref{term2} and \eqref{term3} in (\ref{M1}),  we have
\begin{eqnarray*}
\lefteqn{\frac{1}{2}\partial_t \int_{\RR}  q 
^2\,\le\,\frac{3\,G_0^{\,2}}{2\,\eta }\;e^{T}\;\Big(\Vert p_{0}\Vert_{L^\infty_{y,\sigma}}
+\frac{\sqrt{\alpha\,T}}{\sqrt{\pi}}\,\Big)\vert\partial_y v
 \vert^2}\hskip2truecm\\
 &+&\,\left(\frac{\alpha}{3}+\frac{3}{2\,\eta}\,e^{T}
+\frac{3\,\alpha^{2}}{\eta}\;e^{T}\,\Vert\,\partial_\sigma
p_1\,\Vert^2_{L^2_\sigma}\,\,\right)\;\parallel
q\parallel_{L^{2}_{\sigma}}^2\;,
\end{eqnarray*}
for all $t$ in $[0;T]$ and almost every $y$ in $\Omega$. Applying the  Gronwall lemma,  we obtain
\begin{eqnarray}
\lefteqn{
\int_{\RR}  q ^2(t,y,\sigma)\,d\sigma\le\frac{3\,G_0^{\,2}}{\eta }\,e^{T}\,
\Big(\Vert p_{0}\Vert_{L^\infty_{y,\sigma}}
+\frac{\sqrt{\alpha\,T}}{\sqrt{\pi}}\,\Big)
\int_0^T\vert\partial_y v
 \vert^2\,dt }\nonumber\\
& &\times\,\exp\,\Big(\frac{2\,\alpha\,T}{3}+\frac{3}{\,\eta}\,T\,e^{T}
+\frac{6\,\alpha^{2}}{\eta}\,e^{T}\,\int_0^T\Vert\,\partial_\sigma
p_1\,\Vert^2_{L^2_\sigma}\,dt\Big)\;,
\label{eq:gronq}
\end{eqnarray}
for almost every $(t,y)$ in $\Omega_T$. We now  integrate the above inequality  over $\Omega_T$ 
and we use the bound~\eqref{eq:bdgrad} on $\partial_\sigma p_1$  in $L^\infty_y(L^2_{T,\sigma})$
to deduce \eqref{eq:estiq} with
\begin{eqnarray}
C(T)&=&
\frac{3\,G_0^{\,2}}{\eta}\,T\,e^{T}\,\Big(\Vert p_{0}\Vert_{L^\infty_{y,\sigma}}
+\frac{\sqrt{\alpha\,T}}{\sqrt{\pi}}\,\Big)\;\times\nonumber\\
& & \times\;
\exp\Big(\frac{2\,\alpha\,T}{3}+\frac{3}{\,\eta}\,T\,e^{T}
+\frac{6\,\alpha^{2}}{\eta}\,e^{T}\,C'(T)\Big)
\;,\label{eq:CT}
\end{eqnarray}
with an explicit expression for $C'(T)$ coming from the right-hand side of \eqref{eq:bdgrad}, namely~:
$$C'(T)=\frac{2}{\eta}\,e^{T}\,\left(\Vert
p_0\Vert_{L^\infty_{y,\sigma}}\Big(\frac{1}{2}+T\Big)+
\frac{\alpha}{\sqrt\pi}\,T^{3/2} \right)\;.$$  
\hfill$\diamondsuit$\vskip6pt

\subsection{Global-in-time existence}
\label{sec:longtime} \vskip10pt
Let  us assume that there exists some finite $t^*$ such that the system
admits a solution  $(u^*;p^*)$ on $[0,t^*[$ that ceases to exist after
the time $t^*$ (at least in the appropriate functional spaces prescribed
by our notion of solution). 
\vskip6pt
 Let $0<t_0<t^*$. We consider  the following Cauchy problem starting at time $t_0$~:
 \begin{subequations}
 \label{eq:utilde}
 \begin{EqSystem}
\rho\,\partial_{t}u-\partial^2_{yy} u=\partial_{y}\tau-\rho\,\dot{V}(t)\,y \;;\\
\tau (t,y)=\int_{\RR}\sigma\,p\,d\sigma\;;\\
  u(0,y)=\tilde u_{0}(y):=u^*(t_0, y)\quad \textrm{ on } \Omega\;;\\
    u(t,0)=0\;,\quad u(t,1)=0\;,
\end{EqSystem}
\end{subequations}
which is coupled to 
\begin{subequations}\label{eq:ptilde}
\begin{EqSystem}
\partial_tp\,+\,G_0\,\big(\partial_yu+V(t)\big)\,\partial_\sigma
p-D(p(t,y))\,\partial^2_{\sigma\sigma}p\,
+\un_{\RR\setminus [-1,1]}(\sigma)\,p
=\frac{D(p(t,y))}{\alpha}\,\delta_0(\sigma) \;;
\\ p\,\ge\,0\;;\\
p(0,y,\sigma)=\tilde p_0(y,\sigma):=p^*(t_0, y,\sigma)\;.
  \end{EqSystem}
\end{subequations}
% \begin{remarque}
% Actually the physical situation we are interested in involves zero as the initial condition for $u$ (\textit{i.e.} $u_{0}=0$\/).  It is 
% precisely for the purpose of the argument that we develop now  that we have also  considered general initial conditions $u_{0}$ in $L^2(\Omega)$. \end{remarque}
As a first step we prove that the new initial condition $(\tilde u_0; \tilde p_0)$ satisfies the assumptions in Theorem~\ref{main}. The only point to 
be checked is that  $\tilde p_0$ fulfills~\eqref{strictpos}. Indeed,  $\tilde  u_0$ is in $L^2(\Omega)$ 
and $\tilde p_{0}$ satisfies \eqref{ICfull},   thanks to the bounds in Lemma~\ref{lem:borneL2} which hold true for  
$p^*$ on $[0,t^*)$. Actually we  prove the
more general~:
\begin{lemme}\label{lem:pnondeg} Let  $(u^*;p^*)$ be the solution to the coupled system on $[0;t^*)$  under the assumptions of Theorem~\ref{main} 
for the initial data. Then, for every $0<T<t^*$, we have
\begin{equation}\label{eq:pnondeg}
\inf_{\scriptstyle (t,y)\in\Omega_T\atop \scriptstyle\chi\in\RR}\int_{\vert\sigma+\chi\vert>1}p^*(t,y,\sigma)\,d\sigma\geq
\frac{\eta}{2\,\alpha}\,e^{-T}\;.
\end{equation}
\end{lemme}
\vskip6pt
\noindent\textbf{Proof of Lemma~\ref{lem:pnondeg}: } The proof follows from a comparison principle 
and is inspired from~\cite{CCY}. It is reproduced here for the reader's convenience. 
We denote by $p_{-}$ the solution to the linear equation~:
\begin{equation}\label{eq:systp-}
  \left\{
\begin{array}{rcl}
\partial_{t} p_{-}&=&-\,G_0\left(\partial_y\,u+V(t)\right)\,
  \partial_{\sigma}
p_{-}+D(p^*(t,y))\,\partial_{\sigma\sigma}^{2}
  p_{-}\,-\,p_{-}\;;\\
p_{-}(0,y,\sigma)&=&p_0(y,\sigma)\;.
\end{array}
\right.
\end{equation}
It is well-known that $p_{-}$ is given by
\begin{equation}\label{eq:p-}
p_{-}(t,y,\sigma)= e^{ -t }
\int_{-\infty}^{+\infty} p_{0}(y,\sigma')\,
\varphi_{\sqrt{2\,\int_0^tD(p^*(s,y))\,ds}}(\sigma-\sigma'-\xi(t,y) )\,d\sigma'\;,
\end{equation}
with
\begin{equation*}
\left\{
\begin{array}{rcl}
\varphi_\nu(x)&=&\displaystyle{\frac{1}{\sqrt{2\pi}\;\nu}}\exp\big(-\frac{x^2}{2\,\nu^2}\big)
\quad \mathrm{ if \;}\nu>0\;;\\
\varphi_0&=&\delta_{0} \;,
\end{array}
\right.
\end{equation*}
and 
\begin{equation*}
\xi(t,y)=G_0\,\int_0^t \big(\partial_y\,u(s,y)+V(s)\big)\,ds\;.
\end{equation*}
Since $p_{-}\leq p^*$,  by the maximum principle,  we get the bound from below 
\begin{eqnarray}
\lefteqn{\int_{|\sigma-\chi|>1}p^*(t,y,\sigma)\,d\sigma}\nonumber\\
&\geq &\int_{\vert\sigma-\chi\vert>1} p_{-}(t,y,\sigma)\,d\sigma
\nonumber\\
&\geq&
 e^{-t }\,
\int_\RR p_{0}(y,\sigma')\,\Big( \int_{\vert \sigma-\chi
\vert>1}
\varphi_{\sqrt{2\,\int_0^t (D(p^*(s,y))\,ds}}(\sigma-\sigma'-\xi(t,y))
 \,d\sigma\Big)\,d\sigma'\;,
 \label{eq:strict}
\end{eqnarray}
for every $\chi$ in $\RR$. As in \cite{CCY}, we introduce the interval 
$\displaystyle{K_{\xi,\chi}\,=\,[\,-1-\xi(t,y)+\chi,1-\xi(t,y)+\chi\,]}$.
The function $\sigma\mapsto \varphi_{\sqrt{2\,\int_0^t
D(p^*(s,y))\,ds}}(\sigma-\sigma'-\xi(t,y))$ is a 
Gaussian probability density with mean
$\sigma'+\xi(t,y)$ and
squared width $2\,\int_{0}^{t}\,D(p^*(s,y))\,ds$. 
Therefore, for every
$\displaystyle{\sigma'\in\RR\setminus K_{\xi,\chi}}$,  we have
$$\int_{\vert\sigma -\chi\vert>1}\,\varphi_{\sqrt{2\,\int_0^t
D(p^*(s,y))\,ds}}(\sigma-\sigma'-\xi(t,y))\,d\sigma\ge\frac{1}{2}\;,$$
which implies
\begin{equation*}
\eqref{eq:strict}\geq \frac{1}{2}\,e^{-T}\,\int_{\RR\setminus 
K_{\xi,\chi}}\,p_0(y,\sigma')\,d\sigma'=
\frac{1}{2}\,e^{-T}
 \,\int_{\vert\sigma'-\chi+\xi(t,y)\vert>1} p_0(y,\sigma')\,d\sigma'\;.
\end{equation*} 
And we conclude using~\eqref{strictpos}. 
\hfill$\diamondsuit$
\vskip10pt
\noindent{\bf Completion of the Proof of Theorem~\ref{main}}

\medskip

In view of Lemma~\ref{lem:pnondeg},  we may apply the Banach fixed point theorem as in the proof of Proposition~\ref{prop:ptfixe} and deduce the 
existence of  a unique  solution to the  Cauchy problem \eqref{eq:utilde}--\eqref{eq:ptilde} on the time interval  $[t_0;t_0+\kappa]$ for some 
 small enough $\kappa>0$. Of course, this solution coincides with
 $(u^*;p^*)$ by uniqueness. We now show that
 $\kappa$ may be chosen independently of $t_0$ in  $(0;t^*)$. Therefore
 the solution exists beyond the time $t^*$, which contradicts the
 finiteness of $t^*$. The constant $\kappa$ depends on the Lipschitz
 constant for the mapping ${\cal F}$. Because of  \eqref{eq:lipF2} the Lipschitz constant of the mapping 
$\mathcal{F}_2$ is clearly independent of the initial time $t_0$ and can
be made arbitrarily small using \eqref{eq:lipF2} provided the length of
the time interval is taken small. We thus focus on the Lipschitz
constant of $\mathcal{F}_1$, and now show it is bounded uniformly in
$t^*$. Thus, the condition on $\kappa$ such that $\mathcal{F}=\mathcal{F}_2\circ\mathcal{F}_1$ is a contraction on $[t_0;t_0+\kappa]$
is independent of $t_0$ in $(0;t^*)$, which concludes the proof.

 We revisit carefully the proof  of 
Lemma~\ref{lem:F1}, which is the crucial step for checking the assumptions of the Banach  fixed point theorem on small time interval. 
We go back to the proof of \eqref{eq:estiq}. The only modifications are in  the proofs of  estimates \eqref{term1},  \eqref{term2} and  \eqref{term3} as follows. 
In view of  the uniform estimate given by Lemma~\ref{lem:pnondeg}, the
quantity $D(p^*)$ is bounded from below by $\eta\,\exp(-t^*)\,/2$  and the
$L^\infty$ and $H^1$ norms of $p^*$ in the sense of
\eqref{eq:borneLinfty} and \eqref{eq:bdgrad} are bounded  
uniformly in terms of $t^*$. (In all these bounds $T$ may obviously be bounded by $t^*$.) Therefore the Lipschitz constant of $\mathcal{F}_1$ given by \eqref{cstlip} is bounded uniformly in $t^*$. 
\hfill$\diamondsuit$
\vskip10pt

%%%%%%%%%%%%%%%%%%%%%%%%%%%%%%%%%%%%%%%%%%%%%%%%%%%%%%%%%
\section{The case $\sigma_c=0$}
\label{sec:zero}
%%%%%%%%%%%%%%%%%%%%%%%%%%%%%%%%%%%%%%%%%%%%%%%%%%%%%%%%%%%%%
In the situation examined so far, that is when $\sigma_c>0$, we have
only been able to show the well-posedness of the coupled system when the
macroscopic equation has a \emph{positive} diffusion coefficient
$\mu$. This is a mathematical artefact, apparently related to our
technique of proof. Our aim in the present section, as announced in the
introduction, is to mention that the coupled system is also well-posed
in the particular case when $\mu\geq 0$ \emph{and} $\sigma_c=0$.

We again scale out 
the variables $y$ and $t$, together with the function $u$ as explained
in Appendix, which amounts to taking $T_{0}=L=1$.  
In the present case when $\sigma_c=0$, we note that, for a given $b(t,y)
\in L^2_{\rm loc}(\RR^+,L^2(\Omega))$, the unique solution
of \eqref{eq:syst-p} provided by Proposition~\ref{prop:main1} reads
\begin{eqnarray}
p(t,y,\sigma) & = &  e^{-t} \, \int_\RR p_0(\sigma') \, 
{\cal G}_{\alpha t} \left( \sigma-\sigma'-\chi(t,y) \right) \, d\sigma' 
\label{eq:p-deg} \\
& & + \int_0^t  \, e^{-(t-s)} \, {\cal G}_{\alpha (t-s)}\left(
  \sigma - \chi(t,y) + \chi(s,y) \right) \, ds,  \nonumber
\end{eqnarray}
where $\dps \chi(t,y) = \int_0^t b(s,y) \, ds$ and where $\left( {\cal
    G}_{t} \right)_{t \ge 0}$ denotes the heat kernel 
\begin{eqnarray*}
{\cal G}_t(\sigma)&=&\displaystyle{\frac{1}{\sqrt{4\pi
      t}}}\exp\big(-\frac{\sigma^2}{4\,t}\big)  
\quad \mathrm{ if \;} t>0\;;\\
{\cal G}_0(\sigma) &=&\delta_{0}(\sigma) \;.
\end{eqnarray*}
If we multiply (\ref{eq:p-deg}) by $\sigma$ and integrate over the real
line, we obtain 
\begin{equation}\label{eq:tau-bis}
\partial_t \tau(t,y)  + \tau(t,y)  =  b(t,y) \; .
\end{equation}
For $\sigma_c=0$, the multiscale H\'ebraud--Lequeux model is therefore
equivalent to the so-called Maxwell model~\cite{Owens}, and the coupled
system under consideration reads 
\begin{equation}
  \label{eq:utau}
  \left\{\begin{array}{l}
\displaystyle{\rho\,\partial_t u  =
\mu\, \partial^2_{yy} u + \partial_y \tau - \rho \dot V(t) \, y} \; ; \\
\\
\displaystyle{\partial_t \tau  + \tau  =  G_0\, \partial_y u + G_0 V(t)}
\; ;
  \\ \\
u(\cdot, 0) = u( \cdot, 1) = 0 \; .
  \end{array}
\right.
\end{equation}
The latter is a \emph{linear} system for which it is easy to prove global
existence and uniqueness in convenient functional spaces such that
$\partial_y u \in  L^2_{\rm loc}(\RR^+;L^2(\Omega))$, whatever $\mu\geq 0$.
Global existence and uniqueness for the multiscale
H\'ebraud--Lequeux model immediately follows.

\section*{Appendix: Non-dimensionalized equations}

We first scale the space and time variables in order to work with dimensionless 
constants and with a reduced number of parameters. We introduce  the new dimensionless variables 
\begin{equation*}
t'=\frac{t}{T_{0}}\,,\quad y'=\frac{y}{L}\,,\quad \sigma'=\frac{\sigma}{\sigma_{c}}\;,
\end{equation*} 
and the dimensionless rescaled functions 
\[
U'=\frac{T_{0}}{L}\,U\,,\quad p'=\sigma_{c}\, p\,,\quad \tau'=\frac{\tau}{\sigma_{c}}=\int_\RR \sigma'\,p'\,d\sigma'\;,
\]
together with the corresponding dimensionless parameters
\begin{equation*} 
\rho'=\frac{\rho\,L^2}{\sigma_{c}\,T_{0}{}^2}\,,\quad \alpha'=\frac{\alpha}{\sigma_c^2}\,,
\quad G'_0=\frac{G_0}{\sigma_c}\,, \quad \mu'=\frac{\mu}{T_{0}\,\sigma_{c}}\;.
\end{equation*}
Note that $\rho'$ is actually the so-called Reynolds number. We also define 
\begin{equation*} 
D'(p')=\alpha'\,\int_{|\sigma'|>1}p'\,d\sigma'\;.
\end{equation*}
Then, equations \eqref{eq:u} and \eqref{eq:p} respectively read~:
\be\label{eq:u-res}
\rho\,\partial_{t}U'-\mu\,\partial^2_{yy}U'
=\partial_{y}\tau'\ee
and
\be\label{eq:p-res}
\partial_t p'=-G_0\,\partial_y U '
  \,\partial_\sigma p'+D(p')\,\partial_{\sigma\sigma}^2
p'-\,\un_{\RR \setminus
    [-1,1]}(\sigma)\,p'  +\frac{D(p')}{\alpha}\,\delta_0(\sigma)
\;.
\ee
 Of course the corresponding change of scales and variables are also applied to the initial conditions 
$u_0$ and $p_0$. The new function $V'$ entering in the boundary conditions of $U'$ has to be changed according to 
\ben 
V'= \frac{T_0}{L}\, V\,,
\een
All the primes are omitted in the body of the article in order to lighten
the notation.

\bibliographystyle{plain}

\end{document}